\documentclass[10pt]{amsart}
\usepackage{amsfonts,latexsym}
\usepackage{amsmath}
\usepackage{amscd}
\usepackage{float,amsmath,amssymb,mathrsfs,bm,multirow,graphics}
\usepackage[dvips]{graphicx}
\usepackage[percent]{overpic}
\usepackage[all]{xy}

\usepackage{cite}
\usepackage{enumerate}
\usepackage[percent]{overpic}
\usepackage[all]{xy} 
\usepackage{wrapfig}
  \usepackage{epsfig} 
\usepackage{hyperref}

\newcommand{\nequation}{\setcounter{equation}{0}}
\renewcommand{\theequation}{\mbox{\arabic{section}.\arabic{equation}}}
\newcommand{\R}{\mathbb{R}}

\newcommand{\Z}{\mathbb{Z}}

\newcommand{\s}{{S^1}}

\newcommand{\dx}{{\rm d}x}

\DeclareMathOperator{\supp}{supp}

\newtheorem{theorem}{Theorem}[section]

\newtheorem{definition}[theorem]{Definition}

\newtheorem{remark}[theorem]{Remark}

\input epsf 
\title[On the weakly dissipative Camassa-Holm equation]{On the weakly dissipative Camassa-Holm, Degasperis-Procesi, and Novikov equations}
\author{J. Lenells}
\author{M. Wunsch}
\address{Department of Mathematics, Baylor University, One Bear Place \#97328, Waco, TX
76798, USA}
\address{Departement Mathematik, ETH Z\"urich, R\"amistrasse 101, 8092 Zurich, Switzerland}
\begin{document}

%

\begin{abstract}
We show that the weakly dissipative Camassa-Holm, Degasperis-Procesi, Hunter-Saxton, and Novikov equations can be reduced to their non-dissipative versions by means of an exponentially time-dependent scaling. 
Hence, up to a simple change of variables, the non-dissipative and dissipative versions of these equations are equivalent. 
Similar results hold also for the equations in the so-called $b$-family of equations as well as for the two-component and $\mu$-versions of the above equations.
\end{abstract}



\maketitle

\section{Introduction}
\noindent
The purpose of this note is to point out that the weakly dissipative versions of the Camassa-Holm and several similar equations are equivalent to their non-dissipative counterparts up to a simple change of variables. 

The Camassa-Holm equation
\begin{align}\label{CH}
  u_t-u_{txx}+3uu_x = 2u_xu_{xx}+uu_{xxx}, \qquad x\in\R,\;t>0,
\end{align}
where $u(t,x)$ is a real-valued function, was first obtained in \cite{FF1981} as an abstract bi-Hamiltonian equation and was later derived in the context of water waves \cite{cam-hol, con-lan}. Some of its interesting features are: (a) It is completely integrable \cite{cam-hol, FF1981}. (b) It admits peaked traveling wave solutions (peakons) \cite{cam-hol}. (c) It arises naturally as a geodesic equation in the context of the diffeomorphism group of the circle, see \cite{M1998}. (d) Some solutions exist globally whereas other solutions break in finite time cf. \cite{CE2000}. 

The weakly dissipative Camassa-Holm equation \cite{G2008, N2012, S2011, wu, wu-yin:ch, wu-yin:ch:blowup, YW2010} is obtained by adding to the left-hand side of (\ref{CH}) a term of the form $\lambda (u - u_{xx})$, where $\lambda > 0$ is a parameter:
\begin{align}\label{dissCH}
  v_t-v_{txx}+3vv_x + \lambda (v - v_{xx}) = 2v_xv_{xx}+vv_{xxx}, \qquad x\in\R,\;t>0.
\end{align}
Equation (\ref{dissCH}) has been analyzed in several papers: Local well-posedness of (\ref{dissCH}) in $H^s$, $s > 3/2$, was proved in \cite{wu-yin:ch}.
Existence and uniqueness of global weak solutions of (\ref{dissCH}) was established in \cite{wu}. 
Further global existence and blow-up results were derived in \cite{wu-yin:ch:blowup}. 

However, as we point out in this note, the dissipative and non-dissipative equations (\ref{CH}) and (\ref{dissCH}) are equivalent up to a simple change of variables. More precisely, if $u(t,x)$ and $v(t,x)$ are related by
\begin{align}\label{transformation}
v(t,x)  = e^{-\lambda t} \; u\left(\frac{1 - e^{-\lambda t}}{\lambda}, x\right),
\end{align}
then $u(t,x)$ satisfies (\ref{CH}) if and only if $v(t,x)$ satisfies (\ref{dissCH}).
Using this change of variables, the results available for the non-dissipative equation (\ref{CH}) can be immediately transferred to the dissipative equation (\ref{dissCH}); thus no separate proofs are needed for the weakly dissipative case. 
We will also show that an analogous change of variables exists for many other equations with weakly dissipative versions, such as the Degasperis-Procesi \cite{deg-pro}, Hunter-Saxton \cite{hun-sax}, and Novikov \cite{novikov} equations, the equations in the $b$-family cf. \cite{esc-yin}, 
as well as the two-component \cite{liu-yin, popo} and $\mu$-versions \cite{kh-le-mi} of the above equations. An incomplete list of references that deal with weakly dissipative equations of the type considered in this note includes \cite{hu, guo, gu-la-wa, kohl, niu-zha, N2012, S2011, wu, wu-es-yi, wu-yin:ch, wu-yin:ch:blowup, wu-yin:wd:dp, ya-li-zh, YW2010, zhu-jia}.

It was observed already in \cite{wu-yin:ch:blowup} that the properties of (\ref{dissCH}) are similar to the properties of (\ref{CH}) restricted to a finite time interval. It was also noted that there exist considerable differences between (\ref{CH}) and (\ref{dissCH}) in the long time behavior of solutions; in particular, the global solutions of (\ref{dissCH}) decay to zero as $t\to \infty$. These observations are easy to understand in the light of the correspondence (\ref{transformation}). For example, the decay of the global solutions of (\ref{dissCH}) is an obvious consequence of the exponential prefactor $e^{-\lambda t}$. 

The main results are stated in Section \ref{resultsec}. The proofs are presented in Section \ref{proofsec}. In Section \ref{weaksec} we give an example to illustrate that the correspondence between the dissipative and non-dissipative equations extends also to the case of weak solutions. Finally, in the appendix we motivate the form of the transformation (\ref{transformation}) by deriving explicit solutions of the weakly dissipative Hunter-Saxton system.  

\section{Statement of results}\label{resultsec}\nequation 
\noindent
The weakly dissipative two-component $b$-family of equations \cite{liu-yin} is defined by
\begin{align} \label{wb-sys}
\begin{cases}
m_t + vm_x + b v_x m + \lambda m + \kappa \sigma \sigma_x &= 0,   \\
\sigma_t + (v\sigma)_x + \lambda \sigma &= 0, 
\end{cases}
\qquad x \in \R, \;t > 0,
\end{align}
where $b \in \R$, $\kappa = \pm 1$, and $\lambda \geq 0$ are parameters, $\{v(t,x), \sigma(t,x)\}$ are real-valued functions, and $m = v - v_{xx}$. The system (\ref{wb-sys}) reduces to (\ref{dissCH}) when $b=2$ and $\sigma(t,x) \equiv 0$. 
When $b = 3$ and $\sigma(t,x) \equiv 0$, it reduces to the weakly dissipative Degasperis-Procesi  (DP) equation \cite{guo, gu-la-wa, kohl, wu-es-yi, wu-yin:wd:dp}. If the second component $\sigma(t,x)$ is allowed to be nonzero, the choices $b = 2$ and $b = 3$ in (\ref{wb-sys}) give rise to the weakly dissipative two-component Camassa-Holm and Degasperis-Procesi equations, cf. \cite{hu}. 

Another class of equations is obtained by letting $m = -v_{xx}$ in \eqref{wb-sys}. With this choice, the system \eqref{wb-sys} reduces for $b = 2$ to a weakly dissipative version of the Hunter-Saxton system \cite{len-lec, len-wun,wu-wun,wun:hs,wun:weak}. 

Finally, if we consider spatially periodic solutions and let $m = \mu(v) - v_{xx}$ in (\ref{wb-sys}), where $\mu(v) = \int_\s v \, dx$ denotes the mean of the function $v$ over $S^1 = \R/\Z$, then \eqref{wb-sys} reduces to a weakly dissipative version of the two-component $\mu b$-equation \cite{kh-le-mi}.

The respective non-dissipative equations are obtained by setting $\lambda = 0$ in (\ref{wb-sys}).

\begin{theorem}\label{th1}
Consider the system \eqref{wb-sys} with $m$ defined by either $(i)$ $m = v - v_{xx}$, $(ii)$ $m = -v_{xx}$, or (iii) $m = \mu(v) - v_{xx}$.
Let $(u(t,x), \rho(t,x))$ and $(v(t,x), \sigma(t,x))$ be related by
\begin{align}\label{vurel} 
v(t,x)  = e^{-\lambda t} \; u\left(\frac{1 - e^{-\lambda t}}{\lambda}, x\right), \qquad
\sigma(t,x)  = e^{-\lambda t} \; \rho \left(\frac{1 - e^{-\lambda t}}{\lambda}, x\right). 
\end{align}  
Then $(v,\sigma)$ satisfies the weakly dissipative system \eqref{wb-sys} if and only if $(u,\rho)$ satisfies the respective non-dissipative system. 
\end{theorem}

\begin{remark}\upshape
By choosing $\sigma \equiv \rho \equiv 0$ in (\ref{vurel}), it is clear that the result of Theorem \ref{th1} for the two-component system \eqref{wb-sys} applies also to the corresponding single-component equations.
\end{remark}

\noindent
The weakly dissipative Novikov equation
\begin{align} \label{wd-novi}
 v_t - v_{txx} + 4 v^2 v_x - 3 v v_x v_{xx} - v^2 v_{xxx} + \lambda (v - v_{xx}) = 0, \qquad \lambda > 0, 
\end{align}
was recently studied in \cite{ya-li-zh}. Originally, the Novikov equation
\begin{align} \label{novi}
 u_t - u_{txx} + 4 u^2 u_x - 3 u u_x u_{xx} - u^2 u_{xxx} = 0
\end{align}
was introduced without a dissipative term in \cite{novikov} and subsequently studied in \cite{tiglay:novi,wu-yin:novi}.

\begin{theorem}\label{th2}
Let $u(t,x)$ and $v(t,x)$ be related by
\begin{align}\label{vurel:novi} 
v(t,x)  = e^{-\lambda t} \; u\left(\frac{1 - e^{-2\lambda t}}{2\lambda}, x\right). 
\end{align}  
Then $v$ satisfies the weakly dissipative Novikov equation \eqref{wd-novi} if and only if $u$ satisfies the Novikov equation \eqref{novi}. 
\end{theorem}

\begin{remark}\upshape
The dissipative and non-dissipative solutions in (\ref{vurel}) coincide at time $t = 0$. In other words, the initial data of the two solutions are the same. 
\end{remark}

\begin{remark}\upshape
Under the transformation (\ref{vurel}), dissipative solutions $(v, \sigma)$ on the time interval $[0, T)$ correspond to non-dissipative solutions $(u, \rho)$ on the time interval $[0, S)$ where
\begin{align}\label{ST}
  S = \frac{1 - e^{-\lambda T}}{\lambda} \leq  \frac{1}{\lambda}.
\end{align}
In particular, the dissipative solution $(v, \sigma)$ exists globally if and only if $S \geq 1/\lambda$, where $S$ denotes the maximal existence time of $(u, \rho)$. 
\end{remark}

\section{Proofs}\label{proofsec}\nequation

\noindent
We will present the proof of Theorem \ref{th1}; the proof of Theorem \ref{th2} is similar.

\subsection{Proof of Theorem \ref{th1}}
If not specified otherwise, we will assume that the functions $u, \rho$ (as well as their derivatives) are evaluated at $(\frac{1 - e^{-\lambda t}}{\lambda},x)$, whereas $v, \sigma$ and their derivatives are evaluated at $(t,x)$. 

\textbf{(1)} We first consider the case when $m = v - v_{xx}$. Suppose that $(v,\sigma)$  satisfies (\ref{wb-sys}) with $\lambda > 0$ and that $(u, \rho)$ and $(v, \sigma)$  are related as in (\ref{vurel}). 
Using the notation $n = u - u_{xx}$, the first component of the system \eqref{wb-sys} yields
\begin{align}\nonumber
&\underbrace{-\lambda e^{-\lambda t} \; u + e^{-2\lambda t} \; u_t}_{v_t} 
 + \underbrace{\lambda e^{-\lambda t}\; u_{xx} - e^{-2\lambda t}\; u_{txx}}_{-v_{txx}}
 \\ \label{case1equation}
&+ \underbrace{e^{-2\lambda t} \; u n_x}_{vm_x} + \underbrace{b\; e^{-2\lambda t}\; u_x n}_{bv_x m} +\underbrace{ \lambda \;e^{-\lambda t} \; n}_{ \lambda m} + \underbrace{ \kappa\; e^{-2\lambda t} \; \rho \rho_x}_{ \kappa \sigma \sigma_x}  = 0. 
\end{align}
Multiplying both sides by $e^{2\lambda t}$ and simplifying, we obtain 
$$
n_t + u n_x + b u_x n + \kappa \rho \rho_x = 0, 
$$
which is precisely the first component of \eqref{wb-sys} in the absence of dissipation. 
Similarly, the second component of \eqref{wb-sys} yields
\begin{align*}
\underbrace{-\lambda e^{-\lambda t} \; \rho +  e^{-2\lambda t} \; \rho_t}_{\sigma_t} + \underbrace{e^{-2\lambda t} \; (u \rho)_x}_{(v\sigma)_x} +  \underbrace{\lambda\;e^{-\lambda t}\; \rho}_{\lambda \sigma} = 0. 
\end{align*}
Multiplying both sides by $e^{2\lambda t}$ and simplifying, we obtain 
$$
\rho_t + (u\rho)_x = 0, 
$$
which is the second component of \eqref{wb-sys} in the absence of dissipation.
This proves that $(u, \rho)$ satisfies the non-dissipative system. 
Conversely, if $(u, \rho)$ is a solution of the non-dissipative system, by tracing the above steps backwards, we deduce that $(v, \sigma)$ satisfies  \eqref{wb-sys}  with $\lambda >0$.

\par 
\textbf{(2)} In the case when $m = \mu(v) - v_{xx}$, the analog of equation (\ref{case1equation}) is
\begin{align*}
&\underbrace{- \lambda e^{-\lambda t} \mu(u)
 + e^{-2\lambda t} \mu(u_t)}_{\mu(v_t)} 
 + \underbrace{\lambda e^{-\lambda t}\; u_{xx} - e^{-2\lambda t}\; u_{txx}}_{-v_{txx}} 
\\
&+ \underbrace{e^{-2\lambda t} \; u n_x}_{vm_x} + \underbrace{b\; e^{-2\lambda t}\; u_x n}_{bv_x m} +\underbrace{ \lambda \;e^{-\lambda t} \; n}_{ \lambda m} + \underbrace{ \kappa e^{-2\lambda t} \; \rho \rho_x}_{ \kappa \sigma \sigma_x} = 0,
\end{align*}
where now $n = \mu(u) - u_{xx}$. 
Again, multiplying both sides by $e^{2\lambda t}$ and simplifying, we obtain the first component of \eqref{wb-sys} without dissipation. The proof for the second component follows as in case \textbf{(1)}. 
\par 
\textbf{(3)} The proof when $m = -v_{xx}$ follows by setting $\mu(v) = 0$ in the proof of case \textbf{(2)}.

\section{Weak solutions}\label{weaksec}\nequation
\noindent
The result of Theorem \ref{th1} can be used to transfer results from the non-dissipative to the weakly dissipative setting also in the context of weak solutions. In order to illustrate this point, we here rederive a result which was proved in \cite{wu-yin:ch:blowup} for the weakly dissipative Camassa-Holm equation (\ref{dissCH}) from the analogous result for (\ref{CH}).

Equation (\ref{dissCH}) can be written in weak form as 
\begin{align}\label{weakdissCH}
\begin{cases} 
& v_t + vv_x + \partial_x (1- \partial_x^2)^{-1}\left(v^2 + \frac{1}{2} v_x^2 \right) + \lambda v = 0, 
\\
& v(0,x) = v_0(x),
\end{cases} \qquad x\in\R,\;t>0.
\end{align}

\noindent
The following local well-posedness result was proved in \cite{wu-yin:ch} (see \cite{wu-yin:ch:blowup} for the periodic setting).

\begin{theorem}[Lemma 2.1 of \cite{wu-yin:ch}]
Given $v_0 \in H^r$, $r > 3/2$, there exist a maximal $T = T(v_0) > 0$ and a unique solution $v$ of (\ref{weakdissCH}) such that
\begin{align}\label{vCC1}   
  v = v(\cdot, v_0) \in C([0, T); H^r) \cap C^1([0,T); H^{r-1}).
\end{align}
Moreover, the solution depends continuously on the initial data, i.e. the mapping $v_0 \mapsto v(\cdot, v_0): H^r \to C([0, T); H^r) \cap C^1([0,T); H^{r-1})$ is continuous. 
\end{theorem}

\noindent
Using Theorem \ref{th1}, this result can be derived from the analogous result \cite{R2001} for the Camassa-Holm equation. Indeed, given $v_0 \in H^r$, $r > 3/2$, let $u \in C([0, S); H^r) \cap C^1([0,S); H^{r-1})$ be the solution of (\ref{weakdissCH}) with $\lambda = 0$ and initial data $u(0,x) = v_0(x)$ defined on some maximal time interval $[0, S)$. Define $v(t,x)$ by (\ref{vurel}) for $t \in [0, T)$ where $T = \infty$ if $S \geq 1/\lambda$ and $T = -\ln(1-\lambda S)/\lambda$ if $S < 1/\lambda$. Since $\frac{1 - e^{-\lambda t}}{\lambda}$ is a smooth function of $t \in [0, T)$, it follows from (\ref{vurel}) that $v$ satisfies (\ref{vCC1}). The identity 
$$v_t(t,x) + \lambda v(t,x) = e^{-2\lambda t} u_t\biggl(\frac{1 - e^{-\lambda t}}{\lambda},x\biggr),$$
implies that $v$ satisfies (\ref{weakdissCH}). 
Moreover, since the map $u(\cdot) \mapsto u(\frac{1 - e^{-\lambda \cdot}}{\lambda})$ is continuous
$$C([0, S); H^r) \cap C^1([0,S); H^{r-1}) \to C([0, T); H^r) \cap C^1([0,T); H^{r-1}),$$
the continuous dependence on the initial data for $v$ follows from the analogous property of $u$. 

\begin{remark}\upshape
The consistency of the weak existence results of \cite{wahlen} for the Camassa-Holm equation and of \cite{wu} for its weakly dissipative version, as well as of the results of \cite{es-li-yi:dp} for the Degasperis-Procesi equation and of \cite{gu-la-wa} for its weakly dissipative version, can also be easily established by means of relation \eqref{vurel}. 
\end{remark}


\appendix
\section{The weakly dissipative Hunter-Saxton system}\nequation
\renewcommand{\theequation}{A.\arabic{equation}}
\noindent
The form of the relation (\ref{vurel}) can be motivated by considering the spatially periodic weakly dissipative Hunter-Saxton system:
\begin{align} \label{whs-sys}
\begin{cases}
v_{txx} + v v_{xxx} + 2 v_x v_{xx} + \lambda v_{xx}&= \kappa \sigma \sigma_x, 
\\
\sigma_t + (v \sigma)_x + \lambda \sigma & = 0, 
\\
\qquad v(0,x) &= v_0(x),
\\
\qquad \sigma(0,x) &= \sigma_0(x),
\end{cases}
\qquad x \in S^1, \; t> 0, 
\end{align}
where $\kappa = \pm 1$ and $\lambda > 0$. 
The evolution of the first component in (\ref{whs-sys}) can be rephrased as 
\begin{equation*}
v_{tx} + v v_{xx} + \frac{1}{2} v_x^2 + \lambda v_x = \frac{\kappa}{2} \sigma^2 -2c(t), 
\end{equation*}
where $c(t)$ is determined by periodicity: 
\begin{equation} \label{ct}
c(t) = \frac{1}{4} \int_\s (v_x(t,x)^2 + \kappa \; \sigma(t,x)^2)\;\dx. 
\end{equation}
A straightforward calculation shows that 
\begin{equation*}
e^{2\lambda t} \int_\s \left[ v_x(t,x)^2 + \kappa \sigma(t,x)^2 \right]\;\dx =  \int_\s \left[ v_{0x}(x)^2 + \kappa \sigma_0(x)^2 \right]\;\dx.
\end{equation*}
The Lagrangian flow map $\varphi(t,x)$ is defined by  
\begin{equation} \label{lagrange}
\varphi_t(t,x) = v(t,\varphi(t,x)), \quad \varphi(0,x) = x \in \s. 
\end{equation}
Assuming for simplicity that $c(0) = 1$, slightly tedious calculations reveal that the solution of \eqref{lagrange} is given by  
$$\varphi(t,x) = \int_0^x \biggl\{ \biggl(\cos\tau_\lambda+  \frac{ v_{0x}(x)}{2} \sin\tau_\lambda\biggr)^2 + \kappa \; \frac{\rho_0^2}{4} \sin^2 \tau_\lambda\biggr\} \dx, $$
where we have used the short-hand notation
$$\tau_\lambda := \frac{1 - e^{-\lambda t}}{\lambda}. $$
This leads to the following explicit solution formulas:
\begin{subequations}\label{2HSsolution}
\begin{align} 
v_x(t,\varphi(t,x)) &= e^{-\lambda t} \frac{4 \cos(2 \tau_\lambda) v_{0x} + \sin(2\tau_\lambda) \left( v_{0x}^2 + \kappa \rho_0^2  - 4\right)}{\left[2\cos\tau_\lambda+ v_{0x} \sin \tau_\lambda\right]^2 + \kappa \rho_0^2 \sin^2 \tau_\lambda},
\\[.2cm] 
\rho(t,\varphi(t,x)) &= e^{-\lambda t} \frac{4 \rho_0(x)}{\left[2\cos\tau_\lambda+ v_{0x} \sin \tau_\lambda\right]^2 + \kappa \rho_0^2 \sin^2 \tau_\lambda}. 
\end{align}
\end{subequations}

\noindent
Apart from the prefactor $e^{-\lambda t}$, the solution (\ref{2HSsolution}) coincides exactly with the analogous solution of the non-dissipative system under the correspondence $t \leftrightarrow \tau_\lambda$ (cf. equation (2.6) in \cite{wun:weak}). Thus, defining new functions in terms of the weakly dissipative solution evaluated at $\tau_\lambda$ instead of at $t$ turns these functions into solutions of the original Hunter-Saxton system {\it without} weak dissipation. 
This suggests the relation (\ref{vurel}).

\section*{Acknowledgments}
\noindent
J.L. is grateful for support from the EPSRC.
M.W. acknowledges financial support from the ETH Foundation.

\end{document}